\documentclass{amsart}

\usepackage[utf8]{inputenc}

\usepackage{url}

\usepackage[all]{xy}
\usepackage{pb-diagram, pb-xy}

\usepackage{amssymb}
\usepackage{amsthm}

\newtheorem{thm}{Theorem}[section]
\newtheorem{lem}[thm]{Lemma}
\newtheorem{cor}[thm]{Corollary}

\newtheorem{rem}[thm]{Remark}
\newtheorem{ex}[thm]{Example}

\theoremstyle{plain}
\newtheorem*{gv*}{Assumption $\mathrm{GV}(X)$}

\newcommand{\bbC}{{\mathbb C}}

\newcommand{\bbN}{{\mathbb N}}

\newcommand{\bbZ}{{\mathbb Z}}

\newcommand{\Qbar}{\overline{\mathbb Q}}

\newcommand{\calD}{{\mathcal D}}

\newcommand{\calH}{{\mathcal H}}

\newcommand{\calM}{{\mathcal M}}

\newcommand{\calS}{{\mathcal S}}

\newcommand{\Rep}{{\mathit{Rep}}}

\newcommand{\Dbc}{{{D^{\hspace{0.01em}b}_{\hspace{-0.13em}c} \hspace{-0.05em} }}}
\newcommand{\Perv}{\mathrm{Perv}}

\newcommand{\pH}{{^p \! H}}

\newcommand{\Hom}{{\mathit{Hom}}}

\newcommand{\id}{{\mathit{id}}}

\newcommand{\one}{{\mathbf{1}}}

\newcommand{\overbar}[1]{\mkern 1.5mu\overline{\mkern-1.5mu#1\mkern-3mu}\mkern 3mu}

\newcommand{\bfDbar}{{\overbar{\mathbf{D}}}}

\newcommand{\bfPbar}{{\overbar{\mathbf{P}}}}
\newcommand{\Psibar}{{\overbar{\Psi}}}

\newcommand{\bfPint}{{{\mathbf{P}}_{\mathit{int}}}}
\newcommand{\conv}{\, {*_{\mathit{int}}} \, }

\newcommand{\etabar}{{\mkern 2.5mu \overline{\mkern -2.5mu \eta \mkern -2.5mu} \mkern 2.5mu}}
\newcommand{\sbar}{{\mkern 2mu \overline{\mkern -2mu  s } }}
\newcommand{\Sbar}{{\mkern 2mu \overline{\mkern -2mu S } }}
\newcommand{\Xbar}{{\mkern 2mu \overline{\mkern -2mu X \mkern -2mu} \mkern 2mu}}
\newcommand{\Ybar}{{\mkern 2mu \overline{\mkern -2mu Y \mkern -2mu} \mkern 2mu}}
\newcommand{\Zbar}{{\mkern 2mu \overline{\mkern -2mu Z \mkern -2mu} \mkern 2mu}}
\newcommand{\ibar}{{\overline i}}
\newcommand{\jbar}{{\overline{j}}}

\newcommand{\bfD}{{\mathbf{D}}}

\newcommand{\bfP}{{\mathbf{P}}}

\newcommand{\bfS}{{\mathbf{S}}}
\newcommand{\bfT}{{\mathbf{T}}}

\begin{document}

\title{Perverse sheaves on semiabelian varieties} 
\author{Thomas Kr\"amer}
\address{Mathematisches Institut\\ Ruprecht-Karls-Universit\"at Heidelberg\\ Im Neuenheimer Feld 288, D-69120 Heidelberg, Germany}
\email{tkraemer@mathi.uni-heidelberg.de}

\keywords{Perverse sheaf, convolution product, semiabelian variety, Tannakian category, generic vanishing theorem}
\subjclass[2010]{Primary 14K05; Secondary 14F05, 14F17, 32S60.}

\begin{abstract}
We give a Tannakian description for categories of $l$-adic perverse sheaves on semiabelian varieties which combines a construction of Gabber and Loeser for algebraic tori with a generic vanishing theorem for the cohomology of constructible sheaves on abelian varieties. As an application we explain how the arising Tannaka groups on abelian varieties can be studied in terms of semiabelian degenerations via the functor of nearby cycles.
\end{abstract}

\maketitle

\thispagestyle{empty}

\section{Introduction}

Gabber and Loeser~\cite{GaL} have introduced a Tannakian description for categories of perverse sheaves on algebraic tori over an algebraically closed field $k$ with respect to a certain convolution product. The fibre functor for their Tannakian categories is constructed as the hypercohomology of a generic character twist, using Artin's affine vanishing theorem. For abelian varieties the situation is reversed --- here Artin's theorem no longer applies, but at least over $k=\bbC$ Tannakian arguments can be used to obtain a generic vanishing theorem and a description of perverse sheaves that parallels the one for tori~\cite{KrWVanishing}. The arising Tannaka groups for abelian varieties not only appear in the proof of the vanishing theorem, they also provide a new tool for the study of smooth projective varieties with non-trivial Albanese morphism as illustrated by the case of Brill-Noether theory~\cite{WeBN}, and they can be used to approach classical moduli questions such as the Schottky problem~\cite{KrWSchottky}. So far the most effective way to determine these Tannaka groups has been to consider degenerations of perverse sheaves, and it is desirable to allow also degenerations of the underlying abelian variety into a semiabelian variety. In this note we extend the above constructions to the semiabelian case, and for semiabelian degenerations we show that the nearby cycles functor yields an embedding of the degenerate Tannaka group into the generic one whenever this can possibly hold, see theorem~\ref{thm:nearby}.

\medskip

Before we come to the details, let us give a brief overview over the constructions that follow. Throughout we work over an algebraically closed field $k$ of arbitrary characteristic $p \geq 0$. Let $X$ be a semiabelian variety, i.e.~a commutative group variety which is an extension
$
 1 \rightarrow T \rightarrow X \rightarrow A \rightarrow 1
$
of an abelian variety $A$ by an algebraic torus $T$ over $k$. In what follows we write $\Lambda = \Qbar_l$ for some fixed prime number $l\neq p$ and we denote by
\[
 \bfD \;=\; \bfD(X) \;=\; \Dbc(X, \Lambda) \quad \textnormal{and} \quad  \bfP \;=\; \bfP(X) \;=\; \Perv(X, \Lambda) 
\]
the derived category of bounded constructible complexes of $\Lambda$-sheaves resp.~its full abelian subcategory of perverse sheaves in the sense of~\cite{BBD}. On the derived category the group law~$m: X\times_k X \longrightarrow X$ induces two {\em convolution products} 
\[
 K *_! L \;=\; Rm_! (K\boxtimes L) \quad \textnormal{and} \quad K *_* L \;=\; Rm_* (K\boxtimes L) \quad \textnormal{for} \quad K, L \in \bfD,
\]
but neither of these preserves the abelian subcategory of perverse sheaves. Let us denote by $\bfT \subset \bfD$ the full subcategory of all complexes $K\in \bfD$ which are {\em negligible} in the sense that for all $n\in \bbZ$ the perverse cohomology sheaves $\pH^n(K)$ have Euler characteristic zero. Over the base field $k=\bbC$ we will see that
\begin{enumerate}
\smallskip
\item[(a)] the triangulated quotient category $\bfDbar = \bfD/\bfT$ exists (corollary~\ref{lem:S}), 
 
\smallskip
\item[(b)] both $*_!$ and $*_*$ descend to the same bifunctor $*: \bfDbar \times \bfDbar \rightarrow \bfDbar$ which preserves the essential image
$\bfPbar \subset  \bfDbar$ of $\bfP$
(theorem~\ref{thm:main}), 

\smallskip
\item[(c)] with $*$ as its tensor product, the category $\bfPbar$ is an inductive limit of neutral Tannakian categories (corollary~\ref{cor:tannakian}).
\end{enumerate}

\smallskip
\noindent
Concerning (c) recall from~\cite[th.~2.11]{DM} that a {\em neutral Tannakian category} is a rigid symmetric monoidal $\Lambda$-linear abelian category equivalent to the category $\Rep_\Lambda(G)$ of linear representations of an affine group scheme~$G$ over $\Lambda$. In particular, to any perverse sheaf $P \in \bfP$ we attach in corollary~\ref{cor:tannakian} an affine algebraic group $G=G(P)$ which controls how the convolution powers of the perverse sheaf decompose into irreducible pieces. These groups $G(P)$ are the Tannaka groups we are interested in and whose behaviour under degenerations will be studied in section~\ref{sec:cycles}.

\medskip

In the special case of algebraic tori $X=T$ the above properties (a), (b), (c) have been established by Gabber and Loeser in~\cite[sect.~3]{GaL}. Our arguments are heavily indebted by loc.~cit.~but additionally require the generic vanishing theorem of~\cite{KrWVanishing} for abelian varieties $X=A$. Note that an independent proof of this generic vanishing theorem has been given by Schnell in~\cite{Schn} via the Fourier-Mukai transform for holonomic $\calD$-modules, perhaps closer in spirit to the Mellin transform of Gabber and Loeser but apparently confined to the base field $k=\bbC$. The proof in~\cite{KrWVanishing} also requires the theory of $\calD$-modules, but only at a single step which might generalize to positive characteristic. So below we formulate the generic vanishing theorem for abelian varieties as an axiomatic assumption under which our constructions will work over an algebraically closed field $k$ of arbitrary characteristic.

\medskip

\section{Generic vanishing theorems} \label{sec:gvt}

To formulate a conjectural generalization of the generic vanishing theorem of~\cite{KrWVanishing} to the non-proper case over an algebraically closed field $k$ of characteristic $p\geq 0$ we consider characters of the tame fundamental group. As a base point we take the neutral element $0\in X(k)$. Recall that any semiabelian variety $X$ has a smooth compactification with a normal crossing boundary divisor~\cite[sect.~1.3]{SeQP}.
The tame fundamental group $\pi_1^{\, t}(X,0)$ classifies the finite \'etale coverings of~$X$ that are tamely ramified along each component of such a boundary divisor, see~\cite[exp.~XIII]{SGA1}, \cite{GrothMurre} or~\cite{SchmidtTame}. Thus $\pi_1^{\, t}(X,0)$ is a quotient of the usual \'etale fundamental group $\pi_1(X,0)$, and the two groups are equal if $p=0$ or if $X$ is proper.

\medskip

The group $\Pi(X)$ of continuous characters 
$ \chi:  \pi_1^{\, t}(X,0) \longrightarrow \Lambda^* $
admits a natural decomposition as a direct product
\[ 
 \Pi(X) \;=\; \Pi(X)_{l'} \times \Pi(X)_l ,
\]
where $\Pi(X)_{l'}$ denotes the subgroup of all torsion characters of order prime to~$l$ and where the subgroup~$\Pi(X)_l$ consists of all characters that factor over the maximal pro-$l$-quotient $\pi_1(X,0)_l = \pi_1^{\, t}(X, 0)_l$. This latter quotient is a free \mbox{$\bbZ_l$-module} of finite rank~\cite{BS}, so with the same arguments as in~\cite[sect.~3.2]{GaL} we can identify $\Pi(X)_l$ with the set of~$\Lambda$-valued points of a scheme in a natural way --- though for the same reason as in loc.~cit.~the multiplication of characters does not come from a group scheme structure. Considering 
$\Pi(X)$
as the disjoint union of the infinitely many components $\chi \cdot \Pi(X)_l$ indexed by the characters $\chi \in \Pi(X)_{l'}$, we will say that a statement holds for a {\em generic character} if it holds for all characters in an open subset of $\Pi(X)$ that is dense in each component. 

\medskip

Every character $\chi \in \Pi(X)$ gives rise to a local system $L_\chi$ of rank one on $X$. Our Tannakian constructions are based on the following generic vanishing assumption for the hypercohomology of perverse sheaves, where we denote by $K_\chi = K\otimes_\Lambda L_\chi$ the twist of a sheaf complex~$K\in \bfD(X)$. 

\begin{gv*}
For any perverse sheaf $P\in \bfP(X)$ and generic $\chi \in \Pi(X)$ the forget support morphism
\[ H^\bullet_c(X, P_\chi) \; \longrightarrow \; H^\bullet(X, P_\chi) \] 
is an isomorphism, and in all degrees $i\neq 0$ we have $H^i(X, P_\chi) = H^i_c(X, P_\chi) = 0$.
\end{gv*}

For complex abelian varieties this holds by the vanishing theorem of~\cite{KrWVanishing}. We do not know whether $GV(X)$ is also satisfied for abelian varieties over an algebraically closed field $k$ of characteristic $p=char(k)>0$, but in any case the semiabelian version can be deduced from the abelian one as follows.

\begin{thm} \label{thm:gvt}
If the maximal abelian variety quotient $A = X/T$ of $X$ satisfies the assumption $GV(A)$, then also $GV(X)$ holds.
\end{thm}

{\em Proof.} For any given perverse sheaf $P\in \bfP(X)$, consider with verbatim the same definition as in~\cite[sect.~3.3]{GaL} the Mellin transforms $\calM_!(P)$ and $\calM_*(P)$. These Mellin transforms are objects of the bounded derived category of coherent sheaves on $\Pi(X)$ such that for all $i\in \bbZ$ and any character $\chi \in \Pi(X)$ we have 
%
%\[
 \begin{align*}
 H_c^i(X, P_\chi) \; &\cong \; \calH^i(Li_\chi^*\calM_!(P)), \\
 H^i(X, P_\chi) \; &\cong \; \calH^i(Li_\chi^*\calM_*(P)),
 \end{align*}
%\]
%
where $i_\chi: \{\chi \} \hookrightarrow \Pi(X)$ denotes the embedding of the closed point given by~$\chi$. We also have a morphism $\calM_!(P) \rightarrow \calM_*(P)$ which induces via the above identifications the forget support morphism on cohomology. Since on a Noetherian scheme the support of any coherent sheaf is a closed subset, it follows that the locus of all characters $\chi$ which violate $GV(X)$ forms a closed subset $\calS(P) \subseteq \Pi(X)$.

\medskip

We must show that under the assumption $GV(A)$ the complement of this closed subset $\calS(P)$ intersects every irreducible component of $\Pi(X)$. So we must see that for at least one character $\chi$ in each component the properties in $GV(X)$ hold. To this end consider the exact sequence
\[
 0 \; \longrightarrow \; T \; \stackrel{i}{\longrightarrow} \; X \;
 \stackrel{f}{\longrightarrow} \; A \; \longrightarrow \; 0
\]
which defines our semiabelian variety. Lemma~\ref{lem:torsor} below shows that any component of $\Pi(X)$ contains a character $\chi$ such that the morphism
$
 Rf_!(P_\chi)  \rightarrow  Rf_*(P_\chi)
$
is an isomorphism. Artin's affine vanishing theorem for $f$ then also shows that these two isomorphic direct image complexes are perverse~\cite[lemma~2.4]{GaL}. For all $\varphi\in \Pi(A)$ and $\psi = f^*(\varphi) \in \Pi(X)$ the projection formula says
\[ 
 (Rf_!(P_\chi))_\varphi \;=\; Rf_!(P_{\chi  \psi}) 
 \quad \textnormal{and} \quad 
 (Rf_*(P_\chi))_\varphi \;=\; Rf_*(P_{\chi  \psi}),
\] 
hence we are finished by an application of the vanishing assumption $GV(A)$. \qed

\medskip

To fill in the missing statement in the above proof, we consider our semiabelian variety as a torsor $f: X \rightarrow A$ under $T$ as in~\cite[sect.~III.4]{MiEC}. It will be convenient to forget about semiabelian varieties for a moment. So let $B$ be any variety over $k$, fix a point $b\in B(k)$, and with respect to a torus $T$ consider a torsor
\[
 f: \;\; Y \; \rightarrow \; B \quad \textnormal{with fibre} \quad i: \;\; Y_b \;=\; f^{-1}(b) \; \rightarrow \; Y \, .
\]
Fixing a base point~$y\in Y_b(k)$ we denote by $\Pi(Y_b) = \Pi(T)$ the group of continuous characters of the fundamental group~$\pi_1^t(Y_b, y)$,  and similarly for~$\Pi(Y)$.

\begin{lem} \label{lem:torsor}
In the above setting, for any $P\in \bfP(Y)$ there is a subset $U\subseteq \Pi(Y_b)$ which meets every irreducible component and has the property that the forget support map
$
 Rf_!(P_\chi) \rightarrow Rf_*(P_\chi)
$
is an isomorphism for all $\chi \in \Pi(Y)$ with $i^*(\chi) \in U$.
\end{lem}

{\em Proof.} Any torsor is locally trivial in the \'etale topology~\cite[rem.~4.8(a)]{MiEC}, and the question whether or not $Rf_!(P_\chi) \rightarrow Rf_*(P_\chi)$ is an isomorphism can be checked \'etale-locally on $B$. So via the smooth base change theorem one reduces our claim to the case of the trivial torsor which is treated in~\cite[cor.~2.3.2]{GaL}. 
\qed

\medskip

Returning to our semiabelian variety $X$, for later reference it will be convenient to reformulate theorem~\ref{thm:gvt} in the following equivalent way.

\begin{cor} \label{cor:coh_with_supp}
Let $K\in \bfD(X)$, and suppose the abelian variety $A = X/T$ satisfies the generic vanishing assumption $GV(A)$. Then for all $r\in \bbZ$ and generic $\chi\in \Pi(X)$ the forget support map and perverse truncations induce isomorphisms
\[ H^r_c(X, K_\chi) \; \cong \;  H^r(X, K_\chi) \; \cong \; H^0(X, \pH^r(K)_\chi). \]
\end{cor}

{\em Proof.} For the perverse cohomology groups $P_r = \pH^r(K)$ and $\chi \in \Pi(X)$ we clearly have $(P_r)_\chi = \pH^r(K_\chi)$. Only finitely many $P_r$ are non-zero, so theorem~\ref{thm:gvt} implies that for generic $\chi$ and all $r,s\in \bbZ$ the groups
$H^s_c(X, (P_r)_\chi)$ and $H^s(X, (P_r)_\chi)$
are isomorphic to each other via the forget support morphism, and non-zero at most for $s=0$. Then the spectral sequences
\begin{align} \nonumber
 E_2^{rs} \;=\; H^s(X, (P_r)_\chi) \; \Longrightarrow \; H^{r+s}(X, K_\chi) \\ \nonumber
 E_2^{rs} \;=\; H_c^s(X, (P_r)_\chi) \; \Longrightarrow \; H_c^{r+s}(X, K_\chi) 
\end{align}
degenerate, and our claim follows since the forget support map between the limit terms is induced from the one between the $E_2$-terms.
\qed

\medskip

\section{Some properties of the convolution products} \label{sec:properties-of-convolution}

At this point it will be convenient to gather some properties of the convolution products $*_!$ and $*_*$ for reference. As in~\cite[sect.~2.1]{WeBN} the category $\bfD = \bfD(X)$ is a symmetric monoidal category with respect to each of these products, and in both cases the unit object $\one$ is the rank one skyscraper sheaf supported in the neutral element of the semiabelian variety $X$. For $K\in \bfD$ we consider the Verdier dual~$D(K)$ and define the {\em adjoint dual} to be the pull-back 
$ K^\vee = (-\id)^* D(K) $
under the morphism~$-\id: X\rightarrow X$ that sends an element to its inverse.

\begin{lem} \label{lem:twist}
For $K, M\in \bfD$ and $\chi \in \Pi(X)$ one has the following isomorphisms.\smallskip
\begin{enumerate} 
 \item[\em (a)] Adjoint duality: 
\[
 \qquad \Hom_\bfD(\one, K^\vee *_* M) \;\cong\; \Hom_\bfD(K, M) \;\cong\; \Hom_\bfD(K *_! M^\vee, \one).
\]
\item[\em (b)] Character twists: \vspace{-0.8em}
\[
 \qquad (K_\chi)^\vee \cong (K^\vee)_\chi,
\]
\[
 \qquad (K *_! M)_\chi \; \cong \; K_\chi *_! M_\chi
 \quad \textnormal{\em and} \quad
 (K *_* M)_\chi \; \cong \; K_\chi *_* M_\chi. \vspace{0.4em}
\]
\item[\em (c)] Verdier duality: 
\[
\qquad \qquad
 D(K*_! M) \; \cong \; D(K) *_* D(M), \quad 
 D(K*_* M) \; \cong \; D(K) *_! D(M).
\]
\item[\em (d)] K\"unneth formulae: 
\begin{eqnarray} \nonumber
 \qquad
 H^\bullet_c(X, K) \otimes_\Lambda H^\bullet_c(X, M) &\stackrel{\sim}{\longrightarrow}& H^\bullet_c(X, K*_! M), \\ \nonumber
  H^\bullet(X, K*_* M) &\stackrel{\sim}{\longrightarrow}& H^\bullet(X, K) \otimes_\Lambda H^\bullet(X, M).
  \bigskip
\end{eqnarray}
\end{enumerate}
\end{lem}

{\em Proof.} Part {\em (a)} follows by adjunction as in~\cite[p.~533]{GaL}. The first identity in~{\em (b)} comes from $RHom(L_\chi, \Lambda_X) = L_{\chi^{-1}}$ and $(-\id)^*(L_{\chi^{-1}}) = L_\chi$, the other two follow as in~\cite[prop.~4.1]{KrWVanishing}. Part {\em (c)} follows from the compatibility of Verdier duality with exterior tensor products, and part {\em (d)} is the K\"unneth isomorphism \cite[exp.~XVII.5.4]{SGA4} respectively its Verdier dual. 
\qed

\medskip

\section{The thick subcategory of negligible objects} \label{sec:thick}

In this section we always assume that the abelian variety $A=X/T$ satisfies the assumption $GV(A)$ of section~\ref{sec:gvt}, in which case theorem~\ref{thm:gvt} says that $GV(X)$ holds as well. For sheaf complexes $K\in \bfD = \bfD(X)$ consider the Euler characteristic

\[
 \chi(K) \;=\; \sum_{i\in \bbZ} (-1)^i \dim_\Lambda(H^i(X, K)) \;=\; \sum_{i\in \bbZ} (-1)^i \dim_\Lambda(H^i_c(X, K)),
\]
where the second equality is due to~\cite{La}. We denote by $\bfS=\bfS(X) \subset \bfP=\bfP(X)$ the full subcategory of all perverse sheaves with Euler characteristic zero, and we define~$\bfT=\bfT(X) \subset \bfD(X)$ to be the full subcategory of all sheaf complexes~$K$ such that the perverse cohomology sheaves~$\pH^n(K)$ lie in $\bfS$ for all $n\in \bbZ$.

\begin{lem} \label{lem:twist-invariance}
Let $K\in \bfD$. Then we have $\chi(K) = \chi(K_\varphi)$ for all $\varphi \in \Pi(X)$, and the following conditions are equivalent: \smallskip
\begin{enumerate}
 \item[\em (a)] The complex $K$ lies in the full subcategory $\bfT$. \smallskip
 \item[\em (b)] We have $H^\bullet(X, K_\varphi)=0$ for generic $\varphi \in \Pi(X)$. \smallskip
 \item[\em (c)] We have $H^\bullet_c(X, K_\varphi)=0$ for generic $\varphi \in \Pi(X)$. \smallskip
\end{enumerate}
\end{lem}

{\em Proof.} To prove the invariance of the Euler characteristic under twists we can by d\'evissage assume that $K$ is a constructible sheaf. Then~\cite[cor.~2.9]{IlEuler} says that for any smooth compactification $j: X\hookrightarrow \bar{X}$ the Euler characteristics of the direct images $j_!(K)$ and $j_!(K_\varphi)$ coincide. Here we use that $\varphi$ is tame, see section~2.6 of loc.~cit. 
It then follows that $\chi(K)=\chi(K_\varphi)$, and corollary~\ref{cor:coh_with_supp} implies that the three conditions {\em (a)} -- {\em (c)} are equivalent.  \qed

\medskip

Recall that a full subcategory of an abelian category is called a {\em Serre subcategory} if it is stable under the formation of subquotients and extensions. More generally a full triangulated subcategory of a triangulated category is said to be {\em thick} if for any morphism $f: K \rightarrow L$ which factors over an object of the subcategory and has its cone in the subcategory, both $K$ and $L$ lie in the subcategory as well. 

\begin{cor} \label{lem:S}
The subcategory $\bfS \subset \bfP$ is Serre, and $\bfT\subset \bfD$ is thick.
\end{cor}

{\em Proof.} Consider a short exact sequence $0\rightarrow P \rightarrow Q \rightarrow R \rightarrow 0$ in the abelian category $\bfP$ of perverse sheaves. Theorem~\ref{thm:gvt} says that for generic $\chi\in \Pi(X)$ the hypercohomology of $P_\chi$, $Q_\chi$ and $R_\chi$ is concentrated in degree zero so that after a generic twist the long exact hypercohomology sequence reduces to a short exact sequence 
$0 \rightarrow H^0(X, P_\chi) \rightarrow H^0(X, Q_\chi) \rightarrow H^0(X, R_\chi) \rightarrow 0$.
Hence $Q \in \bfS$ iff~$P, R\in \bfS$ because of the equivalences in lemma~\ref{lem:twist-invariance}. Therefore $\bfS\subset \bfP$ is a Serre subcategory, and $\bfT \subset \bfD$ is then automatically thick by~\cite[prop.~3.6.1(i)]{GaL}.
\qed

\medskip

Localizing the abelian category $\bfP$ at the class of morphisms whose kernel and cokernel lie in $\bfS$ we can thus form the abelian quotient category $\bfPbar = \bfP/\bfS$ in the sense of~\cite[chap.~III]{Ga}. Recall that by definition this quotient category has the same objects as $\bfP$ and that for any objects $P_1, P_2$ the elements of $\Hom_\bfPbar(P_1, P_2)$ can be represented by equivalence classes of diagrams in $\bfP$ of the form
\[
\xymatrix@R=0.5em@C=2em@M=0.5em{
 & Q \ar[dl]_-{f_1} \ar[dr]^-{f_2} & \\
 P_1 & & P_2
}
\]
where the kernel and cokernel of the morphism $f_1$ lie in $\bfS$. Similarly, localizing the triangulated category $\bfD$ at the class of all morphisms whose cone lies in the thick subcategory $\bfT$ we can form the triangulated quotient category $\bfDbar = \bfD / \bfT$ as described in~\cite[sect.~2.1]{Ne}. These two quotient constructions are compatible in the sense that the perverse $t$-structure on~$\bfD$ induces a $t$-structure on $\bfDbar$ whose core is equivalent to $\bfPbar$ in a natural way~\cite[prop.~3.6.1]{GaL}. 

\begin{lem} \label{lem:tensor_ideal} \label{lem:multiplier}
For $K, M\in \bfD$ the following properties hold. \vspace*{0.1em}
\begin{enumerate}
 \item[\em (a)] The cone of the morphism $K*_! M \longrightarrow K*_* M$ lies in $\bfT$. \smallskip
 \item[\em (b)] If $K$ or $M$ lies in $\bfT$, then $K*_!M$ and $K*_*M$ are also in $\bfT$. \smallskip
 \item[\em (c)] If $K, M\in \bfP$, then $\pH^n(K*_!M), \pH^n(K*_*M) \in \bfS$ for all $n\neq 0$. 
\end{enumerate}
\end{lem}

{\em Proof.} {\em (a)} Let $C$ be the cone of the morphism $K*_! M \to K*_* M$. We must show $H^\bullet(X, C_\chi)=0$ for generic $\chi \in \Pi(X)$. For this it suffices to check that on hypercohomology the map
\[
 H^\bullet(X, K_\chi*_!M_\chi) \longrightarrow H^\bullet(X, K_\chi*_*M_\chi)
\]
is an isomorphism for generic $\chi$. By corollary~\ref{cor:coh_with_supp} we can replace the left hand side with the compactly supported hypercohomology, so it remains to see that the forget support map
\[
 f_{K_\chi \boxtimes L_\chi}: \quad H_c^\bullet(X\times_k X, K_\chi \boxtimes M_\chi) \; \longrightarrow \; H^\bullet(X\times_k X, K_\chi \boxtimes M_\chi)
\]
is an isomorphism for generic $\chi$. For simplicity we will suppress the character twist in what follows, replacing $K$ and $L$ by their twists $K_\chi$ and $L_\chi$. Let $j: X\hookrightarrow \Xbar$ be a compactification, and consider the following diagram where all horizontal arrows are forget support morphisms (the vertical arrows will be discussed below). 
\[
\xymatrix@M=0.5em@C=4em@R=1.5em{
 H_c^\bullet(X\times_k X, K\boxtimes L) \ar[r]|-{\; f_{K\boxtimes L} \; }
 & H^\bullet(X\times_k X, K\boxtimes L) \\
 H^\bullet(\Xbar \times_k \Xbar, (j,j)_!(K\boxtimes L)) \ar[r] \ar@{=}[u] \ar@{}[ur]|{(1)}
 &  H^\bullet(\Xbar \times_k \Xbar, R(j,j)_*(K\boxtimes L)) \ar@{=}[u] \\
 H^\bullet(\Xbar \times_k \Xbar, j_!(K)\boxtimes j_!(L)) \ar[r] \ar[u] \ar@{}[ur]|{(2)}
 &  H^\bullet(\Xbar \times_k \Xbar, Rj_*(K)\boxtimes Rj_*(L)) \ar[u] \\
 H^\bullet(\Xbar, j_!(K)) \otimes H^\bullet(\Xbar, j_!(L)) \ar[r] \ar[u] \ar@{}[ur]|{(3)}
 & H^\bullet(\Xbar, Rj_*(K)) \otimes H^\bullet(\Xbar, Rj_*(L)) \ar[u] \\
 H_c^\bullet(X, K) \otimes H_c^\bullet(X, L) \ar[r]|-{\; f_K \otimes f_L \;} \ar@{=}[u] \ar@{}[ur]|{(4)}
 & H^\bullet(X, K) \otimes H^\bullet(X, L) \ar@{=}[u] 
}
\]
The horizontal arrow $f_K \otimes f_L$ is the tensor product of the forget support morphisms for $K$ and $L$, hence by corollary~\ref{cor:coh_with_supp} it is an isomorphism for generic $\chi$. On the other hand, the horizontal arrow $f_{K\boxtimes L}$ is the forget support morphism we are interested in. So we will be finished if we can show that the above diagram commutes and that all the vertical arrows (to be defined yet) are isomorphisms.

\medskip

The squares (1) and (4) commute by definition of $f_{K\boxtimes L}$ and $f_K \otimes f_L$. The vertical arrows in (3) are the K\"unneth isomorphisms for the proper morphism $\Xbar \to \mathit{Spec}(k)$, and the commutativity of this square follows from the naturality of the K\"unneth isomorphism. Finally, the square (2) is induced by the square
\[
\xymatrix@M=0.5em@C=2em{
 (j,j)_!(K\boxtimes L) \ar[r] & R(j,j)_*(K\boxtimes L) \\
 j_!(K) \boxtimes j_!(L) \ar[r] \ar[u]_-{\mathit{ad}_!^{-1}}  
 & Rj_*(K) \boxtimes Rj_*(L) \ar[u]_-{\mathit{ad}_*}
}
\]
where $\mathit{ad}_!$ and $\mathit{ad}_*$ are the natural morphisms which correspond via adjunction to the identity of
\[
 K\boxtimes L \;=\; (j,j)^!(j_!(K)\boxtimes j_!(L)) \;=\; (j,j)^*(Rj_*(K)\boxtimes Rj_*(L)).
\]
Note that $\mathit{ad}_!$ is an isomorphism and that over $U=X\times_k X \hookrightarrow \Xbar\times_k \Xbar$ all the morphisms in the above diagram restrict to the identity. The diagram commutes since by adjunction there is only one morphism
$
 j_!(K) \boxtimes j_!(L) \rightarrow  R(j,j)_*(K\boxtimes L)
$
which restricts over $U$ to the identity. By the same argument the diagram
\[
\xymatrix@C=5em@M=0.5em@R=1.5em{
 Rj_*(K) \boxtimes Rj_*(L) \ar[r]^-{\mathit{ad}_*} \ar@{=}[d] & R(j,j)_*(K\boxtimes L) \ar@{=}[d] \\
 D(j_!(D(K)) \boxtimes j_!(D(L))) \ar[r]^-{D(\mathit{ad}_!)} & D((j,j)_!(D(K)\boxtimes D(L)))
}
\]
commutes. Since the lower row is an isomorphism, it follows that $\mathit{ad}_*$ must be an isomorphism as well, and this finishes the proof of part {\em (a)} of the lemma.

\medskip

{\em (b)} For $K\in \bfT$ and generic $\chi\in \Pi(X)$ we have $H_c^\bullet(X, K_\chi)=0$ by lemma~\ref{lem:twist-invariance}, and the K\"unneth formula in lemma~\ref{lem:twist} then implies
\[
 H^\bullet_c(X, (K*_!M)_\chi) \;=\; H^\bullet_c(X, K_\chi) \otimes_\Lambda H^\bullet_c(X, M_\chi) \;=\; 0
\]
as well. So $K*_! M$ lies in $\bfT$ as required. The statement for $K*_*M$ follows in the same way or can be deduced via Verdier duality.

\medskip

{\em (c)} For generic $\chi \in \Pi(X)$ theorem~\ref{thm:gvt} says that $H^\bullet_c(X, K_\chi) \cong H^\bullet(X, K_\chi)$ and that this hypercohomology is concentrated in degree zero. The same holds with $M$ in place of $K$. So lemma~\ref{lem:twist} shows that $H^\bullet_c(X, K*_! M)$ and $H^\bullet(X, K*_*M)$ are concentrated in degree zero for generic $\chi$, and we can apply corollary~\ref{cor:coh_with_supp}.
\qed 

\medskip

\section{Tannakian categories} \label{sec:tannakian_semiabelian}

In this section we assume as above that for the semiabelian variety $X$ and its quotient $A=X/T$ the generic vanishing assumption $GV(A)$ and hence also $GV(X)$ are satisfied. The results of the previous section then imply

\begin{thm} \label{thm:main}
The two convolution products $*_!$ and $*_*$ both descend to the same well-defined bifunctor
\[
 *: \; \bfDbar \times \bfDbar  \; \longrightarrow \; \bfDbar
\]
which satisfies $\bfPbar * \bfPbar \subset \bfPbar$ and with respect to which both $\bfDbar$ and $\bfPbar$ are symmetric monoidal \mbox{$\Lambda$-linear} triangulated categories.
\end{thm}

{\em Proof.} By lemma~\ref{lem:tensor_ideal}~{\em (b)} both $*_!$ and $*_*$ descend to a bifunctor on $\bfDbar$, and part~{\em (a)} of the lemma shows that these two bifunctors coincide. As in~\cite[sect.~2.1]{WeBN} it follows that $(\bfDbar, *)$ is a symmetric monoidal category, and part {\em (c)} of lemma~\ref{lem:tensor_ideal} furthermore shows that we have $\bfPbar * \bfPbar \subseteq \bfPbar$. \qed

\medskip

It is sometimes convenient to lift the quotient category $\bfPbar$ to a category of true perverse sheaves. To obtain such a lift, consider the full subcategory $\bfPint \subset \bfP$ of all perverse sheaves without subquotients in~$\bfS$. Then as in~\cite[sect.~3.7]{GaL} the quotient functor $\bfP \rightarrow \bfPbar$ restricts to an equivalence of categories between $\bfPint$ and $\bfPbar$, and via this equivalence the convolution product $*$ induces a bifunctor
\[
 \conv: \quad \bfPint \times \bfPint \; \longrightarrow \; \bfPint
\]
with respect to which $\bfPint$ is a $\Lambda$-linear symmetric monoidal category equivalent to $\bfPbar$. The unit object $\one$ of $\bfPint$ is the rank one skyscraper sheaf supported in the neutral element of the group variety $X$. 
As an application we show that the category $\bfPbar$ is rigid, i.e.~that we have a notion of duality which involves evaluation and coevaluation morphisms with the usual properties~\cite[def.~1.7]{DM}.

\begin{thm} \label{thm:rigid}
The symmetric monoidal abelian category
$\bfPbar$ is rigid.
\end{thm}

{\em Proof.} Let $P\in \bfPint$ be a perverse sheaf without constituents in $\bfS$. Via the adjunction property in lemma~\ref{lem:twist} the identity morphism $\id_P: P\longrightarrow P$ defines two morphisms $\one \longrightarrow P^\vee *_* P$ and $P*_! P^\vee \longrightarrow \one$ in the derived category $\bfD$. Via the quotient functor these two morphisms induce morphisms in the full subcategory $\bfPbar$ by lemma~\ref{lem:multiplier}. Denote by
\[
 coev: \; \one \longrightarrow P^\vee \conv P \qquad \textnormal{and} \qquad ev: \; P \conv P^\vee \longrightarrow \one
\]
the corresponding morphisms in the lift $\bfPint$. We must show that for all $P \in \bfPint$ the composite morphism
\[
\gamma: \;
\xymatrix@C=5.5em@M=0.5em{
 P = P \conv \one  \ar[r]^{\id \conv coev} & P \conv P^\vee \conv P \ar[r]^{ev\conv \id} & \one \conv P = P 
}
\]
and its counterpart $(\id \conv ev) \circ (coev \conv \id): \; P^\vee \longrightarrow P^\vee$ are the identity. Since the argument is the same in both cases, we will only deal with $\gamma$ in what follows.

\medskip

In order to show $\gamma - \id_P = 0$, we can by the axiom $GV(X)$ assume that for all subquotients $Q$ of the perverse sheaf $P$ the hypercohomology $H^\bullet(X, Q)$ vanishes in all non-zero degrees. Then $H^\bullet(X, -)$ behaves like an exact functor on all short exact sequences which only involve subquotients of $P$. After a suitable character twist we can furthermore assume that 
\[ H^\bullet(X, P^\vee *_! P) \;=\; H^\bullet(X, P^\vee \conv P) \;=\; H^\bullet(X, P^\vee *_* P) \]
and that the forget support morphism for these hypercohomology groups is an isomorphism. Then for $H=H^\bullet(X, P)$ and $H^\vee = \Hom_\Lambda(H, \Lambda)$ the diagrams given in \cite[appendix A.5.4]{GaL} show that the morphism
\[
\xymatrix@C=1em{
\Lambda \;=\; H^\bullet(X, \one) \ar[rr]^-{coev_*} && H^\bullet(X, P^\vee \conv P) \;=\; H^\bullet(X, P^\vee *_* P) \;=\; H^\vee \otimes_\Lambda H
}
\]
is the coevaluation in the category of vector spaces, and dually for $ev$. Since the category of vector spaces of finite dimension over $\Lambda$ is rigid, it follows that on hypercohomology $\gamma - \id_P$ induces the zero morphism. Then $Q=P/\ker(\gamma - \id_P)$ is acyclic in the sense that $H^\bullet(X, Q)=0$. By definition of $\bfPint$ it follows that $Q=0$ and hence $\gamma - \id_P = 0$ as required.
\qed

\medskip

For $P\in \bfPbar$ let us now denote by $\langle \, P \, \rangle \subset \bfPbar$ the full subcategory of all objects which are isomorphic to subquotients of convolution powers of $P\oplus P^\vee$. Subcategories of this form are said to be {\em finitely tensor generated}. By construction they inherit from~$\bfPbar$ the structure of a rigid symmetric monoidal $\Lambda$-linear abelian category, and we claim that they are neutral Tannakian:

\begin{cor} \label{cor:tannakian}
For any $P\in \bfPbar$ there is an affine algebraic group~$G=G(P)$ over~$\Lambda$, unique up to non-canonical isomorphism, such that we have an equivalence 
\[
 \omega: \; \; \langle \, P \, \rangle \; \stackrel{\sim}{\longrightarrow} \; \Rep_\Lambda(G)
\]
where $\Rep_\Lambda(G)$ denotes the rigid symmetric monoidal $\Lambda$-linear abelian category of finite-dimensional algebraic representations of $G$. 
\end{cor}

{\em Proof.} Consider $P$ as an object of $\bfPint \subset \bfP$. By theorem~\ref{thm:gvt} there exists a character $\chi \in \Pi(X)$ such that all constituents of all convolution powers 
\[ (P\oplus P^\vee)_\chi \conv \cdots \conv (P\oplus P^\vee)_\chi
\]
have their hypercohomology concentrated in degree zero (note that these are only countably many conditions and that $\Pi(X)_l$ cannot be covered by countably many Zariski-closed subsets). For such a character $\chi$ the functor
$ Q \mapsto H^0(X, Q_\chi) $
is a fibre functor from $\langle \, P \, \rangle$ to the category of finite-dimensional vector spaces, and our corollary follows from theorem~\ref{thm:rigid} via the Tannakian formalism~\cite[th.~2.11]{DM}. 
\qed

\medskip

The above fibre functor depends on the chosen character twist, and this is the reason why we have restricted ourselves to finitely tensor generated subcategories and why the Tannaka group is only determined up to isomorphism. Alternatively one could as in~\cite[th.~3.7.5]{GaL} take the fibre functor which to a perverse sheaf assigns the generic fibre of its Mellin transform, but this fibre functor is defined over a large function field, and the descent to $\Lambda$ again involves non-canonical choices.

\begin{rem} With notations as above, \smallskip
\begin{enumerate}
 \item[\em (a)] if~$P\in \bfP(X)$ is semisimple, then $G(P)$ is a reductive algebraic group, \smallskip
 \item[\em (b)] if $X$ is a complex abelian variety, then $G(P)$ is the group defined in~\cite{KrWVanishing}.
\end{enumerate}
\end{rem}

{\em Proof.} In {\em (a)} the group $G(P)$ acts faithfully on the representation $V=\omega(P)$ which in the case at hand is semisimple. Since the unipotent radical of an algebraic group acts trivially on every irreducible representation, it follows that the unipotent radical of $G(P)$ must be trivial, and we are done. Claim~{\em (b)} follows easily from the universal property of the localization $\bfDbar$. \qed

\medskip

\section{Nearby cycles} \label{sec:cycles}

To describe the behaviour of the above Tannaka groups under degenerations, we work over the spectrum $S$ of a strictly Henselian discrete valuation ring with closed point $s$. We consider perverse sheaves with coefficients in $\Lambda = \Qbar_l$ for a prime $l$ that is invertible on $S$. Fixing a geometric point $\etabar$ over the generic point $\eta\in S$ we denote by $\Sbar$ the normalization of $S$ in the residue field of $\etabar$.

\medskip

Consider a semiabelian scheme $X \to S$, i.e.~a smooth commutative group scheme over $S$ whose fibres $X_s$ and $X_\eta$ are semiabelian varieties.  Throughout this section we always assume that the generic fibre $X_\eta$ is an abelian variety. Putting~$\Xbar  = X \times_S \Sbar$ we have a  commutative diagram
\[
\xymatrix@M=0.4em@R=0.4em@C=1.5em{
 & X_\sbar \ar[rr]^{\quad \ibar} \ar@{-}[d] \ar@{=}[dl] && \Xbar \ar@{-}[d] \ar[dl] && X_\etabar \ar[ll]_{\quad \jbar} \ar[dd] \ar[dl] \\
 X_s \ar[dd] \ar[rr]^-{\qquad i} & \ar[d] & X \ar[dd]  & \ar[d] & X_\eta \ar[dd] \ar[ll]_-{\qquad j} & \\
 & \sbar \ar@{=}[dl] \ar@{-}[r] & \ar[r] & \Sbar \ar[dl] & \ar[l] &  \etabar  \ar@{-}[l] \ar[dl] \\
 s  \ar[rr] && S &&  \eta  \ar[ll]
}
\]
in which all squares are Cartesian. Degenerations of constructible sheaves can in this setting be studied via the functor of nearby cycles~\cite[exp.~XIII-XIV]{SGA7}
\[
 \Psi: \quad \bfD(X_\etabar) \, \longrightarrow \, \bfD(X_s), \quad \Psi(K) \,=\, \ibar^* (R\jbar_* (K)).
\]
By~\cite[th.~4.2 and cor.~4.5]{Il} this functor commutes with Verdier duality and restricts to an exact functor between the abelian categories $\bfP(X_\etabar)$ and $\bfP(X_s)$ of perverse sheaves. However, in the non-proper case the functor of nearby cycles in general does not preserve the Euler characteristic.

\begin{ex} \label{ex:euler}
For $n\in \bbN$ the kernel $X[n] \hookrightarrow X$ of the multiplication by $n$ is a quasi-finite group scheme over $S$. It decomposes as a disjoint union
$
 X[n] = Y \, \amalg \, Z
$
where $Y$ is finite over $S$ and where $Z \cap X_s = \varnothing$. In general $Z \neq \varnothing$ but a support argument shows that the skyscraper sheaf $P=\Lambda_{Z_\etabar} \in \bfP(X_\etabar)$ satisfies $\Psi(P)=0$.
\end{ex}

For the rest of this section we always assume that for all geometric points $t$ of~$S$ the fibre $X_t$ satisfies the generic vanishing assumption $GV(X_t)$ of section~\ref{sec:gvt}. We then consider the quotient categories 
\[ \bfDbar{(X_t)} = \bfD{(X_t)} / \bfT{(X_t)} \quad \textnormal{and} \quad \bfPbar{(X_t)} = \bfP{(X_t)} / \bfS{(X_t)} \]
as defined in section~\ref{sec:tannakian_semiabelian}. Note that the Euler characteristic is well-defined on objects of these quotient categories.

\begin{lem} \label{lem:upper_estimate}
The functor $\Psi$ descends to a functor $\Psibar: \bfPbar(X_\etabar) \to \bfPbar(X_s)$. For all objects $P, Q\in \bfPbar(X_\etabar)$ we have 
\begin{align}
 0 \; \leq \; \chi(\Psibar(P)) \; & \leq \; \chi(P), \label{eq:chi_estimate} \tag{$a$} \\
 0 \; \leq \; \chi ( \Psibar(P) * \Psibar(Q) ) \; & \leq \; \chi ( \Psibar (P * Q) ), \tag{$b$} \label{eq:conv_estimate}
\end{align}
and $\Psibar(P) * \Psibar(Q)$ is a direct summand of $\Psibar(P*Q)$ inside the category $\bfPbar(X_s)$.
\end{lem}

{\em Proof.} We begin with some general remarks. Let $f: Y\to Z$ be a homomorphism of semiabelian $S$-schemes. By abuse of notation we again write $f$ for any base change of it. For the base changes $\Ybar = Y\times_S \Sbar$ and $\Zbar = Z\times_S \Sbar$ we have a commutative diagram
\[
\xymatrix@R=1.8em@C=2.5em@M=0.5em{
 Y_{s} \ar[r]^-\ibar \ar[d] & \Ybar \ar[d] & Y_{\etabar} \ar[l]_-\jbar \ar[d] \\
 Z_{s} \ar[r]^-\ibar & \Zbar & Z_{\etabar} \ar[l]_-\jbar
}
\]
over $\Sbar$. For $K\in \bfD(\Ybar)$ this gives a commutative diagram
\[
\xymatrix@R=1.8em@C=3em@M=0.5em{
 \ibar^* (Rf_! (K)) \ar[r]^-{bc_!} \ar[d]_\alpha & Rf_! (\ibar^* (K)) \ar[d]^\beta \\
 \ibar^* (Rf_* (K)) \ar[r]^{bc_*} & Rf_* (\ibar^* (K))
} 
\]
where $\alpha$ and $\beta$ denote the forget support morphisms and where $bc_!$ and~$bc_*$ are the base change morphisms. Recall how these latter are defined: The direct image of the adjunction morphism $K\rightarrow i_* (i^*(K))$ yields two morphisms
%
%\[
\begin{align*}
 Rf_!(K) \; \longrightarrow \;\; & Rf_! (\ibar_* (\ibar^*(K)))  \; = \;  \ibar_* (Rf_!(\ibar^*(K))), \\
 Rf_*(K) \; \longrightarrow \;\; & Rf_* (\ibar_* (\ibar^*(K)))  \, =  \; \ibar_* (Rf_*(\ibar^*(K))),
\end{align*}
%\]
%
and again by adjunction these give rise to the base change morphisms $bc_!$ and~$bc_*$ from above. We also remark~\cite[exp.~XVII, th.~5.2.6]{SGA4} that as a result of the proper base change theorem $bc_!$ is always an isomorphism. Putting $K=R\jbar_*(L)$ for a complex $L\in \bfD(Y_\etabar)$, we obtain a factorization
\[
\xymatrix{
 Rf_! (\Psi_{Y}(L)) \ar[rr]^-\beta \ar[dr]_-{\alpha \circ (bc_!)^{-1} \; \; } & & Rf_* (\Psi_{Y}(L)) \\
 & \Psi_{Z}(Rf_*(L)) \ar[ur]_-{bc_*} &
}
\]
where $\Psi_Y$ and $\Psi_Z$ denote the nearby cycles for $Y\to S$ resp.~$Z\to S$. We apply these remarks in the following two situations.

\medskip

First we take $f: X\to S$ to be the structure morphism with $Y=X$, $Z=S$ and~$L=P$. The above diagram then more explicitly says that the forget support morphism $\beta: H_c^\bullet(X_s, \Psi(P)) \rightarrow H^\bullet(X_s, \Psi(P))$ factors over $H^\bullet(X_\etabar, P)$. 
Let $\chi$ be a character of the tame fundamental group of $\Xbar$, and $\chi_t\in \Pi(X_t)$ its pull-back to the fibre $X_t$ over a geometric point $t$ of $S$. In the appendix in section~\ref{sec:sp_characters} we will see that
\[
 \Psi(P_{\chi_\etabar}) \;=\; (\Psi(P))_{\chi_s}
\]
and that $\chi$ can be chosen in such a way that the characters $\chi_s$ and $\chi_\etabar$ are both generic. By the generic vanishing axiom $GV(X_s)$ we can therefore assume that the forget support morphism~$\beta$ is an isomorphism, in which case the above factorization shows that $H^\bullet(X_s, \Psi(P))$ is a direct summand of~$H^\bullet(X_\etabar, P)$. Furthermore, by the generic vanishing axiom $GV(X_\etabar)$ we can also assume that all the occuring hypercohomology groups are concentrated in cohomology degree zero. It follows that $0 \leq \chi(\Psi(P)) \leq  \chi(P)$, in particular the nearby cycles functor $\Psi$ sends the Serre subcategory $\bfS(X_\etabar)\subset \bfP(X_\etabar)$ into the Serre subcategory $\bfS(X_s)\subset \bfP(X_s)$ and hence induces a functor
\[
 \Psibar: \; \bfPbar(X_\etabar) \longrightarrow \bfPbar(X_s)
\]
between the quotient categories. The property $(a)$ of the functor $\Psibar$ is inherited from the corresponding estimate for $\Psi$ shown above.

\medskip

Secondly we take $f=m$ to be the group law, with $Y=X\times_k X$, $Z=X$ and~$L=P\boxtimes Q$. In $\bfDbar(X_s)$ we can identify $\Psi_Z(Rf_*(L))$ with $\Psibar(P*Q)$. Since by~\cite[th.~4.7]{Il} the exterior tensor product $\boxtimes$ commutes with nearby cycles, we can also identify $Rf_*(\Psi_Y(L))$ with $\Psibar(P)*\Psibar(Q)$. Then the factorization of $\beta$ from above more explicitly says that the nearby cycles $\Psibar(P*Q)$ admit $\Psibar(P)*\Psibar(Q)$ as a direct summand in~$\bfPbar(X_s)$. Hence claim $(b)$ follows.
\qed

\medskip

In general we cannot expect $\Psibar$ to be compatible with convolution products since it does not preserve the Euler characteristic (see example~\ref{ex:euler}). Let us say that a perverse sheaf~$P\in \bfP(X_\etabar)$ is {\em admissible} if 
\[ \chi(\Psibar(P)) \;=\; \chi(P). \] 
For an abelian scheme $X\to S$ of course every perverse sheaf is admissible since the nearby cycles are compatible with proper morphisms.

\begin{lem} \label{lem:admissible_tensor_category}
The admissible objects form a rigid symmetric monoidal full abelian subcategory 
\[ \bfPbar(X_\etabar)^{ad} \;\subset \; \bfPbar(X_\etabar) \]
which is stable under the formation of extensions and subquotients.
\end{lem}

{\em Proof.} Any exact sequence $0\rightarrow P_1 \rightarrow P_2 \rightarrow P_3 \rightarrow 0$ in $\bfPbar(X_\etabar)$ is mapped under the nearby cycles functor to a short exact sequence in the category $\bfPbar(X_s)$. Then in particular
\[
 \chi(P_2) \;=\; \chi(P_1) + \chi(P_3)
 \quad \textnormal{and} \quad
 \chi(\Psibar(P_2)) \;=\; \chi(\Psibar(P_1)) + \chi(\Psibar(P_3))
\]
since the Euler characteristic is additive in short exact sequences. But on the other hand $\chi(P_i)\geq \chi(\Psibar(P_i))$ by lemma~\ref{lem:upper_estimate} {\em (a)}. Hence $P_2$ is admissible iff $P_1$ and~$P_3$ are both admissible. So the category of admissible objects is stable under extensions and subquotients. It remains to show that for admissible $P_1, P_2 \in \bfPbar(X_\etabar)$ also the convolution product $P_1 * P_2$ is admissible. This follows from
\begin{align} \nonumber
 \chi(P_1 * P_2) & = \chi(P_1) \cdot \chi(P_2) & \textnormal{(part $(d)$ of lemma 3)} \\ \nonumber
 & = \chi(\Psibar(P_1)) \cdot \chi(\Psibar(P_2)) & \textnormal{(admissibility of $P_1$, $P_2$)} \\ \nonumber
 & = \chi(\Psibar(P_1)*\Psibar(P_2)) & \textnormal{(part $(d)$ of lemma 3)} \\ \nonumber
 & \leq \chi(\Psibar(P_1 * P_2)) & \textnormal{(part $(b)$ of lemma 12)} \\ \nonumber
 & \leq \chi(P_1*P_2) & \textnormal{(part $(a)$ of lemma 12)}
\end{align}
which implies that the last two estimates must in fact be equalities. \qed

\medskip

For the compatibility of the nearby cycles with convolution products it now follows that the situation is as good as one could possibly hope for. 

\begin{thm} \label{thm:nearby}
On the rigid symmetric monoidal abelian subcategory of admissible objects, the nearby cycles define a tensor functor ACU in the sense of~\cite{Rivano}
\[
 \Psibar: \quad \bfPbar(X_{\etabar})^{ad} \; \longrightarrow \; \bfPbar(X_s).
\]
In particular, for $P\in \bfPbar_{\etabar}^{ad}$ we have a closed immersion $G(\Psibar(P)) \hookrightarrow G(P)$.
\end{thm}

{\em Proof.} For $P_1, P_2 \in \bfPbar(X_\etabar)^{ad}$ the last statement in lemma~\ref{lem:upper_estimate} {\em (b)} shows that we have in $\bfPbar(X_s)$ a split monomorphism
$\Psibar(P_1)* \Psibar(P_2)  \hookrightarrow  \Psibar(P_1 * P_2)$.
But we know from the proof of lemma~\ref{lem:admissible_tensor_category} that the source and target of this monomorphism have the same Euler characteristic. So the monomorphism is in fact an isomorphism and hence $\Psibar$ is a tensor functor on the category of admissible objects. For the statement about the Tannaka groups recall that every object of $\langle \Psibar(P) \rangle$ is a subquotient of the image under~$\Psibar$ of some object in~$\langle P \rangle$. So our claim follows from the Tannakian formalism~\cite[prop.~2.21b)]{DM}.\qed

\medskip

\section{Appendix: Specialization of characters} \label{sec:sp_characters}

Let $X\rightarrow S$ be a semiabelian scheme as in the previous section. Since $l$ is assumed to be invertible on $S$, we know that
\[
 \pi_1(X_\etabar,0)_l \;=\; \pi_1^{\, t}(X_\etabar,0)_l 
 \quad \textnormal{and} \quad
 \pi_1(X_s,0)_l \;=\; \pi_1^{\, t}(X_s,0)_l
\]
are free $\bbZ_l$-modules of finite rank~\cite{BS}. However, the existence of abelian varieties with semiabelian reduction shows that in general $\pi_1(X_s, 0)_l$ may have smaller rank than $\pi_1(X_\etabar, 0)$. For the proof of lemma~\ref{lem:upper_estimate} we need to justify why in passing from the generic to the special fibre, we retain enough characters to apply the generic vanishing assumption of section~\ref{sec:gvt} on both fibres. Let~$\Pi(\Xbar)_l$ denote the group of continuous characters of $\pi_1(\Xbar, x)_l$ for any geometric point~$x$ in~$\Xbar$. The passage to a different point $x$ corresponds to an inner automorphism of the fundamental group and is not seen on the level of characters. So we have well-defined restriction homomorphisms
%
%\[
\begin{align*}
 \ibar^*: \quad & \Pi(\Xbar)_l \;\longrightarrow \; \Pi(X_s)_l, \quad \chi \mapsto \chi_s, \\
 \jbar^*: \quad & \Pi(\Xbar)_l \;\longrightarrow \; \Pi(X_\etabar)_l, \quad  \chi \mapsto \chi_\etabar. 
 \end{align*}
%\]
%
From the theory of N\'eron models one deduces that $\ibar^*$ is surjective and that $\jbar^*$ is an isomorphism, see~\cite[sect.~3.6]{KrDiss}. This in particular implies that if a statement holds for a generic character on the fibres $X_s$ resp.~$X_\etabar$, then we can find a global character~$\chi \in \Pi(\Xbar)_l$ such that the statement holds for both $\chi_s$ and $\chi_\etabar$. This is precisely what we need for the proof of lemma~\ref{lem:upper_estimate}. The character $\chi$ defines a local system $L_\chi$ on $\Xbar$, and for all $K\in \bfD(X_\etabar)$ the projection formula implies
\medskip
\[
 \Psi(K_{\chi_\etabar}) 
 \;=\; \ibar^*(R\jbar_*(K\otimes_\Lambda \jbar^* (L_\chi)))
 \;=\; \ibar^*(R\jbar_*(K)\otimes_\Lambda L_\chi)
 \;=\; \Psi(K)_{\chi_s}
 \medskip
\]
so that we can apply the generic vanishing axiom from section~\ref{sec:gvt} simultaneously on both fibres $X_\etabar$ and $X_s$ as required.

\bigskip

{\em Acknowledgements.} This note is based on chapter 3 of my Ph.D. thesis, and I would like to thank my advisor R. Weissauer for his continuous interest in my work and for many inspiring mathematical discussions.

\medskip

\bibliographystyle{amsplain}
\bibliography{BibliographySemiabelian}

\end{document}